\documentclass[12pt,twoside]{amsart}
\usepackage{amssymb,}
\usepackage{latexsym,exscale}
\usepackage[english,francais]{babel}
\usepackage{times}
\usepackage{url}
\input utile.def
\newtheorem{theorem}{Theorem}[section]

\newtheorem{e-proposition}[theorem]{Proposition}

\newtheorem{e-definition}[theorem]{Definition\rm}
\newtheorem{remark}{\it Remark\/}

\newtheorem{theoreme}{Th\'eor\`eme}[section]

\newtheorem{remarque}{\it Remarque}

\newcommand{\pa}{\partial}
\def\RR{{\mathbb{R}}}
\begin{document}
\selectlanguage{francais}
\title{%
R\'egularit\'e du rayon hyperbolique }
\author{%
Satyanad Kichenassamy~$^{\text{a}}$ }
\address{%
\begin{itemize}\labelsep=2mm\leftskip=-5mm
\item[$^{\text{a}}$]
Laboratoire de Math\'ematiques (UMR 6056), CNRS \&\ Universit\'e
de Reims Champagne-Ardenne, Moulin de la Housse, B.P. 1039, 51687
Reims Cedex 2 \\ Courriel: {\tt
satyanad.kichenassamy@univ-reims.fr}
\end{itemize}
}
\maketitle

\begin{quote}{\textbf{R\'esum\'e :}
Soit $\Omega\subset\RR^2$ un domaine born\'e de classe
$C^{2+\alpha}$, $0<\alpha<1$. On montre que si $u$ est la solution
maximale de $\Delta u = 4\exp(2u)$, qui tend vers $+\infty$ si
$(x,y)\to\partial\Omega$, alors le rayon hyperbolique $v=\exp(-u)$
est de classe $C^{2+\alpha}$ jusqu'au bord. La d\'emonstration
repose sur de nouvelles estimations de Schauder pour des
\'equations fuchsiennes elliptiques. }
\end{quote}
\selectlanguage{english}
\title{R\'egularit\'e du rayon hyperbolique/Smoothness of hyperbolic radius}
\begin{quote}{\textbf{Abstract:}
Let $\Omega\subset\RR^2$ be a bounded domain of class
$C^{2+\alpha}$, $0<\alpha<1$. We show that if $u$ is the maximal
solution of $\Delta u = 4\exp(2u)$, which tends to $+\infty$ as
$(x,y)\to\partial\Omega$, then the hyperbolic radius $v=\exp(-u)$
is of class $C^{2+\alpha}$ up to the boundary. The proof relies on
new Schauder estimates for Fuchsian elliptic equations.
}\end{quote}


\textbf{Note presented by Ha\"\i m Brezis,} received and accepted for publication on October 26, 2003.

\textbf{Appeared in:} \textit{Comptes Rendus. Math\'ematique}, Volume 338 (2004) no.~1, pp.~13-18. \textbf{DOI} : 10.1016/j.crma.2003.10.037. \url{https://comptes-rendus.academie-sciences.fr/mathematique/articles/10.1016/j.crma.2003.10.037/}

\begin{center}
\textit{Version anglaise abr\'eg\'ee}
\end{center}

\section{Introduction}

Let $\Omega\subset\RR^2$ be a bounded domain of class
$C^{2+\alpha}$, with $0<\alpha<1$. Let $d(x,y)$ denote the
distance of $(x,y)$ to the boundary. Let $u$ be the maximal
solution of
\begin{equation}
-\Delta u +4e^{2u}=0.
\end{equation}
It satisfies $u(x,y)\to \infty$ as $(x,y)\to\pa\Omega$.

The hyperbolic radius $v$ is defined by $v=\exp(-u)$; for its
properties and applications \emph{see} [2]. We prove:
\begin{theorem}\label{th:e1}
$v$ is of class $C^{2+\alpha}$ up to the boundary.
\end{theorem}
\begin{remark} Theorem 1.1 was conjectured
in [2, p.~204]. It is proved in [5, th.~2.4] that $v$ is of class
$C^{2+\beta}$ for some $\beta>0$, provided that $\Omega$ is
(convex and) of class $C^{4+\alpha}$. This enables one to justify
the formal asymptotics $v=2d-d^2(\kappa+o(1))$, where $\kappa$ is
the curvature of $\pa\Omega$, from which it follows, if $\Omega$
is strictly convex, that $u$ is convex near the boundary. Theorem
1.1 is optimal since $\kappa$ is of class $C^\alpha$.
\end{remark}
\begin{remark} An asymptotic expansion of $u$ to second order, for
smooth domains, is given in [3, p.~32].
\end{remark}
\begin{remark}
There is an extensive literature on the issue of boundary blow-up;
see [1--3,7,9--14] and their references for further details.
\end{remark}

The proof rests on the study of the properties of the renormalized
unknown $w$ defined by
\[ w(x,y)=\frac{v- 2d}{d^2}.\]
This new unknown solves the Fuchsian equation
\begin{equation}\label{eq:e-fuchs}
Lw+2\Delta d = M_w(w),
\end{equation}
where
\[ L=d^2\Delta +2d\nabla d\cdot\nabla-2
    =\mathrm{div}(d^2\nabla)-2,\]
and $M_w$ is a linear operator with $w$-dependent coefficients:
for any $f$, we let
\[ M_w(f) = \frac{d^2}{2+dw}\left[2f\nabla w\cdot\nabla d
+ d\nabla w\cdot\nabla f \right] -2df\Delta d.
\]

A synopsis of the proof follows. It uses auxiliary results proved
in sections 3, 4 and 5.

\section{Synopsis of proof}

Let $\Omega'\subset\Omega$ be a thin $C^{2+\alpha}$ domain, on
which $d$ is less than a parameter $\delta\leq 1/2$, and $d$ is of
class $C^{2+\alpha}$. We work throughout on $\Omega'$.

The first step is to establish, using a comparison argument and
interior regularity, that
\begin{theorem}\label{th:e2}
$w$ and $d^2\nabla w$ are bounded near $\pa\Omega$.
\end{theorem}
Theorem 3.1 applied to $L-M_w$ then leads to
\begin{theorem}\label{th:e3}
If $\delta$ is small, then $dw$ and $d^2\nabla w$ are of class
$C^\alpha(\overline{\Omega'})$, and $d\nabla w$ is bounded near
$\pa\Omega$.
\end{theorem}
Next, one finds $w_0$ which is smooth enough, and is such that
$\tilde{w}=w-w_0$ has controlled boundary behavior:
\begin{theorem} \label{th:e4}
If $\delta$ is small, there is a function $w_0$ such that $w_0$,
$d\nabla w_0$, and $d^2\nabla^2 w_0$ belong to
$C^\alpha(\overline{\Omega'})$, and
\[ Lw_0+2\Delta d=0\text{ on }\Omega'.
\]
\end{theorem}
A second comparison argument (\emph{see} section 5), then yields
\begin{theorem}\label{th:e5} There is a constant $\gamma$ such
that
\[ |\tilde{w}|\leq \gamma d\ln (1/d)
\text{ on }\Omega'.
\]
\end{theorem}
Applying first theorem 3.2, then theorem 3.3, to $L$, we find
\begin{theorem}\label{th:e6}
If $\delta$ is small, $d^2\tilde{w}$ is of class
$C^{2+\alpha}(\overline{\Omega'})$.
\end{theorem}
Theorem 1.1 follows.

\section{Scaled Schauder estimates}

An operator $A$ is said to be \emph{of type (I)} (on a given
domain) if it can be written
\[ A=\pa_i(d^2a^{ij}\pa_{j})+db^i\pa_i+c, \]
with $(a^{ij})$ uniformly elliptic and of class $C^\alpha$, and
$b^i$, $c$ bounded.

An operator is said to be \emph{of type (II)} if it can be written
\[ A=d^2a^{ij}\pa_{ij}+db^i\pa_i+c, \]
with $(a^{ij})$ uniformly elliptic and $a^{ij}$, $b^i$, $c$ of
class $C^\alpha$.
\begin{theorem}\label{th:e7}
If $Ag=f$, where $f$ and $g$ are bounded and $A$ is of type (I) on
$\Omega'$, then $d\nabla g$ is bounded, and $dg$ and $d^2\nabla g$
belong to $C^\alpha(\Omega'\cap\pa\Omega)$.
\end{theorem}
\begin{theorem}\label{th:e8}
If $Ag=df$, where $f$ and $g$ are bounded, $g=O(d^\alpha)$, and
$A$ is of type (I) on $\Omega'$, then $g\in
C^{\alpha}(\Omega'\cap\pa\Omega)$, and $dg\in
C^{1+\alpha}(\Omega'\cap\pa\Omega)$
\end{theorem}
\begin{theorem}\label{th:e9}
If $Ag=df$, where $f\in C^\alpha(\overline{\Omega'})$,
$g=O(d^\alpha)$, and $A$ is of type (II) on $\Omega'$, then $d^2
g\in C^{2+\alpha}(\Omega'\cap\pa\Omega)$.
\end{theorem}
These results are variants of those in [6, chap.~6].

\section{Regularity up to the boundary for a model Fuchsian equation}

One can localize the problem, using a partition of unity, to the
neighborhood of a point $P$ on the boundary. Performing a rigid
motion if necessary, we may assume that $\nabla d= (1,0)$ at $P$.
One then performs the change of coordinates $(x,y)\mapsto(T,Y)$,
where $T=d(x,y)$ and $Y=y$. In this coordinate system, $L$ takes
the form $L=L_0+L_1$, where
\[ L_0 = (D+2)(D-1)+T^2\pa^2_Y,  \]
where $D=T\pa_T$. We work on $\{0\leq T\leq \theta, |Y|\leq
\theta\}$, $\theta$ small.

Let $k$ be $2\theta$-periodic in $Y$, and let
\[ \tilde{k}=\int_1^\infty
k(T\sigma,Y)\frac{d\sigma}{\sigma^2},
\]
which solves $(D-1)\tilde{k}=-k$. Next, find $h(x,y)$ by solving
\[\Delta h+\tilde{k}=0,\]
with $h=0$ for $T=0$, $h_T=0$ for $T=\theta$, and periodic
conditions in $Y$. Finally, let
\[w_1=T^{-2}[(D-1)h].\]
This function satisfies $L_0w_1=k$ and has the required smoothness
properties. One then treats equation $(L_0+L_1)w_0=k$ by a
perturbation argument.

\section{Construction of sub- and super-solutions}

Direct computation shows that $-\ln[2d+d^2(w_0+ A d\ln d)]$ is a
super-solution of (1) if $A$ is large and positive, resp.\ a
sub-solution if $A$ is large and negative. One then uses an
argument from [4] (\emph{see also} [12]) which enables one to
conclude $u$ lies between these two functions. Theorem 2.4
follows.

\section{Perspectives}

The methods of this work may be generalized to produce systematic
constructions of comparison functions and expansions, of high
order, along the lines of our earlier work on the hyperbolic case
(\emph{see e.~g.} [8]). The systematic introduction of Fuchsian
operators enables one to construct comparison functions, justify
expansions and determine the optimal regularity of solutions.

\par\medskip\centerline{\rule{2cm}{0.2mm}}\medskip
\setcounter{section}{0}
\setcounter{equation}{0}
\selectlanguage{francais}

\begin{center}\textit{Texte principal (en fran\c cais).}\end{center}

\section{Introduction}

Soit $\Omega\subset\RR^2$ un domaine born\'e de classe
$C^{2+\alpha}$, avec $0<\alpha<1$. On d\'esigne par $d(x,y)$ la
distance de $(x,y)$ au bord. Soit $u$ la solution maximale de
\begin{equation}\label{eq:liouville}
-\Delta u +4e^{2u}=0.
\end{equation}
Elle satisfait $u(x,y)\to \infty$ quand $(x,y)\to\pa\Omega$.

Le rayon hyperbolique $v$ est d\'efini par $v=\exp(-u)$; pour ses
propri\'et\'es et applications \emph{voir} [2]. On montre ici:
\begin{theoreme}\label{th:1}
$v$ est de classe $C^{2+\alpha}$ jusqu'au bord.
\end{theoreme}
\begin{remarque} Le th\'eor\`eme 1.1 a \'et\'e conjectur\'e
dans [2, p.~204]. Dans [5, th.~2.4], on d\'emontre que $v$ est de
classe $C^{2+\beta}$ pour un certain $\beta>0$, si $\Omega$ est
(convexe et) de classe $C^{4+\alpha}$. Ceci permet de justifier le
d\'eveloppement formel $v=2d-d^2(\kappa+o(1))$, o\`u $\kappa$
d\'esigne la courbure du bord, d'o\`u l'on tire, si $\Omega$ est
strictement convexe, que $u$ est convexe pr\`es du bord. Le
th\'eor\`eme 1.1 est optimal car $\kappa$ est de classe
$C^\alpha$.
\end{remarque}
\begin{remarque} Un d\'eveloppement asymptotique de $u$ \`a
l'ordre deux pour des domaines tr\`es r\'eguliers est donn\'e dans
[3, p.~32].
\end{remarque}
\begin{remarque}
Le probl\`eme de l'explosion au bord a \'et\'e \'etudi\'e dans de
nombreux travaux; on pourra consulter en particulier
[1--3,7,9--14] et leurs r\'ef\'erences pour une vue plus
compl\`ete.
\end{remarque}

La d\'emonstration repose sur l'\'etude de l'inconnue
renormalis\'ee $w$ d\'efinie par
\[ w(x,y)=\frac{v-2d}{d^2}.\]
Cette nouvelle inconnue satisfait l'\'equation fuchsienne
\begin{equation}\label{eq:fuchs}
Lw+2\Delta d = M_w(w),
\end{equation}
o\`u
\[ L=d^2\Delta +2d\nabla d\cdot\nabla-2
    =\mathrm{div}(d^2\nabla)-2,\]
et $M_w$ est un op\'erateur lin\'eaire \`a coefficients
d\'ependant de $w$: pour toute $f$, on pose
\[ M_w(f) = \frac{d^2}{2+dw}\left[2f\nabla w\cdot\nabla d
+ d\nabla w\cdot\nabla f \right]-2df\Delta d.
\]
On donne un r\'esum\'e de la preuve, qui utilise des r\'esultats
auxiliaires d\'emontr\'es dans les sections 3, 4, et 5.

\section{R\'esum\'e de la preuve}

On travaille sur un domaine $\Omega'\subset\Omega$, sur lequel $d$
n'exc\`ede pas $\delta\leq 1/2$, que l'on prendra assez petit dans
la suite, et sur lequel $d$ est est de classe $C^{2+\alpha}$.

La premi\`ere \'etape consiste \`a d\'emontrer, par un argument de
comparaison, que
\begin{theoreme}\label{th:2}
$w$ et $d^2\nabla w$ sont born\'ees sur $\Omega'$.
\end{theoreme}
Le th\'eor\`eme 3.1 appliqu\'e \`a $L-M_w$, fournit alors:
\begin{theoreme}\label{th:3}
Si $\delta$ est petit, $dw$ et $d^2\nabla w$ sont dans
$C^\alpha(\overline{\Omega'})$, et $d\nabla w$ est born\'ee sur
$\Omega'$.
\end{theoreme}
Il faut ensuite trouver $w_0$ assez r\'eguli\`ere, telle que
$\tilde{w}=w-w_0$ soit assez plate au bord (sections 4 et 5) :
\begin{theoreme} \label{th:4}
Si $\delta$ est assez petit, il existe une fonction $w_0$ telle
que $w$, $d\nabla w$, et $d^2\nabla^2 w$ sont dans
$C^\alpha(\overline{\Omega'})$ et v\'erifie
\[ Lw_0+2\Delta d=0
\]
sur $\Omega'$.
\end{theoreme}
Un second argument de comparaison (\emph{voir} section 5) fournit
alors :
\begin{theoreme}\label{th:5} Il existe une constante $\gamma$ telle que
\[ |\tilde{w}|\leq \gamma d\ln (1/d)
\]
sur $\Omega'$.
\end{theoreme}
Appliquant le th\'eor\`eme 3.2, puis le th\'eor\`eme 3.3 \`a $L$,
il vient
\begin{theoreme}\label{th:6}
Si $\delta$ est petit, $d^2\tilde{w}$ est de classe
$C^{2+\alpha}(\overline{\Omega'})$.
\end{theoreme}
Le th\'eor\`eme 1.1 en r\'esulte.

\section{Estimations int\'erieures et changements d'\'echelle}

Le th\'eor\`eme 2.1 se d\'emontre par un argument de comparaison.
Les th\'eor\`emes 2.2 et 2.5 r\'esultent de deux r\'esultats
g\'en\'eraux sur les estimations de Schauder pour des op\'erateurs
fuchsiens lin\'eaires, appliqu\'es \`a $L$ ou \`a $L-M_w$, sur
$\Omega'$. On dira qu'un op\'erateur $A$ est \emph{du type (I)}
(sur un ouvert donn\'e) s'il a la forme
\[
   A=\pa_i(d^2a^{ij}\pa_{j})+db^i\pa_i+c,
\]
avec $(a^{ij})$ uniform\'ement elliptique et de classe $C^\alpha$,
et $b^i$, $c$ born\'ees.

Un op\'erateur sera dit \emph{du type (II)} s'il s'\'ecrit
\[
  A=d^2a^{ij}\pa_{ij}+db^i\pa_i+c,
\]
avec $(a^{ij})$ uniform\'ement elliptique et $a^{ij}$, $b^i$, $c$
de classe $C^\alpha$.

On montre alors, par des arguments de changements d'\'echelle \`a
partir des estimations int\'erieures, que
\begin{theoreme}\label{th:7}
Si $Ag=f$, o\`u $f$ et $g$ sont born\'ees, et $A$ est de type (I)
sur $\Omega'$, alors $d\nabla g$ est born\'ee sur $\Omega'$, et
$dg$ et $d^2\nabla g$ sont de classe
$C^\alpha(\overline{\Omega'})$.
\end{theoreme}
\begin{theoreme}\label{th:8}
Si $Ag=df$, o\`u $f$ et $g$ sont born\'ees, $g=O(d^\alpha)$, et
$A$ est de type (I) sur $\Omega'$, alors $g\in
C^{\alpha}(\overline{\Omega'})$ et $dg\in
C^{1+\alpha}(\overline{\Omega'})$
\end{theoreme}
\begin{theoreme}\label{th:9}
Si $Ag=df$, o\`u $f\in C^\alpha(\overline{\Omega'})$,
$g=O(d^\alpha)$, et $A$ de type (II) sur $\Omega'$, alors $d^2 g$
est de classe $C^{2+\alpha}(\overline{\Omega'})$.
\end{theoreme}
Ces r\'esultats sont des variantes de ceux de [6, ch.~6]; on
pourrait sans doute affaiblir encore les hypoth\`eses sur les
coefficients comme dans [6, th.~6.2].

\section{Etude d'une \'equation fuchsienne mod\`ele}

On se ram\`ene, \emph{via} partition de l'unit\'e, \`a un
probl\`eme local, pr\`es d'un point $P$ du bord. On peut, par un
d\'eplacement, supposer qu'en ce point $\nabla d=(1,0)$. On
introduit alors le changement de coordonn\'ees
$(x,y)\mapsto(T,Y)$, o\`u $T=d(x,y)$ et $Y=y$. On trouve alors que
dans ces coordonn\'ees, $L=L_0+L_1$, o\`u
\[ L_0 = (D+2)(D-1)+T^2\pa^2_Y,  \]
avec $D=T\pa_T$. On travaille sur $\{0\leq T\leq \theta, |Y|\leq
\theta\}$, $\theta$ petit.

On suppose $k$ p\'eriodique en $Y$ et l'on consid\`ere
\[ \tilde{k}=\int_1^\infty
k(T\sigma,Y)\frac{d\sigma}{\sigma^2},
\]
qui v\'erifie $(D-1)\tilde{k}=-k$. On obtient ensuite $h(x,y)$
telle que
\[\Delta h+\tilde{k}=0,
\]
avec $h=0$ pour $T=0$, $h_T=0$ pour $T=\theta$, et $h$
p\'eriodique en $Y$. On v\'erifie enfin que
\[w_1=T^{-2}[(D-1)h]
\]
v\'erifie $L_0w_1=k$ et poss\`ede les propri\'et\'es de
r\'egularit\'e exig\'ees au th\'eor\`eme 2.3. On r\'esout ensuite
l'\'equation $(L_0+L_1)w_0=k$ par un argument de perturbation.
Ceci d\'emontre le th\'eor\`eme 2.3.

\section{Construction de sur- et sous-solutions}

On montre, par un calcul direct, que $-\ln[2d+d^2(w_0+ A d\ln d)]$
est une sur-solution de (1) pr\`es du bord si $A$ est positif et
assez grand, et une sous-solution si $A$ est n\'egatif et assez
grand. On utilise alors un argument d\^u \`a [4] (\emph{voir}
\'egalement [12]) qui permet de conclure que $u$ est comprise
entre ces deux fonctions. Le th\'eor\`eme 2.4 en r\'esulte.

\section{Perspectives}

Les m\'ethodes utilis\'ees dans ce travail se g\'en\'eralisent
naturellement dans deux directions. D'une part, les techniques
g\'en\'erales que nous avons utilis\'ees dans le cas hyperbolique
fournissent des m\'ethodes nouvelles de construction de fonctions
de comparaison et de d\'eveloppements asymptotiques d'ordre
\'elev\'e, et pour des non lin\'earit\'es assez g\'en\'erales
(\emph{voir par exemple} [8]). D'autre part, l'introduction
syst\'ematique d'op\'erateurs fuchsiens permet \'egalement de
justifier ces d\'eveloppements et de d\'eterminer la
r\'egularit\'e optimale des solutions.


%
\end{document}